# AFFINELY INVARIANT MATCHING METHODS WITH DISCRIMINANT MIXTURES OF PROPORTIONAL ELLIPSOIDALLY SYMMETRIC DISTRIBUTIONS

By Donald B. Rubin and Elizabeth A. Stuart[1]

*Harvard University and Johns Hopkins University*

In observational studies designed to estimate the effects of interventions or exposures, such as cigarette smoking, it is desirable to try to control background differences between the treated group (e.g., current smokers) and the control group (e.g., never smokers) on covariates $X$ (e.g., age, education). Matched sampling attempts to effect this control by selecting subsets of the treated and control groups with similar distributions of such covariates. This paper examines the consequences of matching using affinely invariant methods when the covariate distributions are "discriminant mixtures of proportional ellipsoidally symmetric" (DMPES) distributions, a class herein defined, which generalizes the ellipsoidal symmetry class of Rubin and Thomas [*Ann. Statist.* **20** (1992) 1079–1093]. The resulting generalized results help indicate why earlier results hold quite well even when the simple assumption of ellipsoidal symmetry is not met [e.g., *Biometrics* **52** (1996) 249–264]. Extensions to conditionally affinely invariant matching with conditionally DMPES distributions are also discussed.

**1. Background.** The goal in many applied projects is to estimate the causal effect of a treatment (e.g., cigarette smoking) from nonrandomized data by comparing outcomes (e.g., lung cancer rates) in treated (e.g., current smokers) and control (e.g., never smokers) groups, after adjusting for covariate differences (e.g., age, education) between the groups. A common method is to form matched subsamples of the treated and control groups such that the distributions of covariates $X$ are more similar in the matched samples than in the original groups. The use of matched sampling has been

Received July 2005; revised September 2005.
[1]Supported in part by a National Science Foundation Graduate Research Fellowship.

receiving more and more attention in fields such as statistics (e.g., [11, 15]), economics (e.g., [4, 7, 10, 21]), political science (e.g., [9]), sociology (e.g., [20]) and medicine (e.g., [1]) as a class of methods for controlling bias in such observational studies. Here we provide theoretical guidance for choosing matching methods that reduce bias in the matched groups, as well as guidance on the amount of bias reduction that can be achieved with fixed distributions and fixed sample sizes.

We begin with random samples from the treated and control groups of fixed sizes $N_t$ and $N_c$, respectively, with $X$ measured in both samples. Matching chooses subsamples of fixed sizes $N_{mt}$ and $N_{mc}$ from the original groups on which to measure the outcome variables, as well as possibly measure additional covariates. Throughout, we use the subscripts $t$ and $c$ to indicate quantities in the original random samples from the treated and control groups, and the subscripts $mt$ and $mc$ to indicate the corresponding quantities in the matched treated and control groups.

We restrict attention to a particular but general class of matching methods, those that are affinely invariant. In practice, many matching methods are affinely invariant in the sense that the same matched samples will be obtained after any full-rank affine transformation of $X$. For example, the same matches will be obtained if people's heights are measured in inches or centimeters, or if their temperatures are measured in degrees Fahrenheit or degrees Kelvin. Formally, let $\mathcal{X}_t$ and $\mathcal{X}_c$ be data matrices (units by variables). A matching method is a mapping from $(\mathcal{X}_t, \mathcal{X}_c)$ to a pair of sets of indices $(T, C)$ representing the units chosen in the matched samples. An affinely invariant matching method results in the same output $(T, C)$ after any (full-rank) affine transformation $A$ of the $X$:

$$(\mathcal{X}_t, \mathcal{X}_c) \to (T, C) \quad \text{implies} \quad (A(\mathcal{X}_t), A(\mathcal{X}_c)) \to (T, C).$$

Affinely invariant matching methods include Mahalanobis metric, discriminant or propensity score matching. Non-affinely invariant methods include methods where one coordinate of $X$ is treated differently from the others or where nonlinear estimators of the discriminant (or other metric) are used, as discussed by Rubin and Thomas [16].

Theoretical results in papers by Rubin and Thomas [16, 17] describe the effects of affinely invariant matching on bias reduction, as well as on variance, in the matched treated and matched control groups, when $X$ has ellipsoidally symmetric distributions (e.g., the normal distribution or the multivariate $t$) in the treated and control groups, with proportional covariances. Rubin and Thomas [18] used these theoretical results to obtain a series of approximations for the bias and variance reduction possible in a particular matching setting using true and estimated propensity scores, with no subsampling of the treated sample and normal distributions. They then examined the performance of these approximations by simulation with ellipsoidal nonnormal distributions and found that the approximations based on the normal distribution held remarkably well, even for a $t$-distribution with 5 degrees of freedom. They also explored the performance of the approximations with real data from a study of prenatal hormone exposure, with 15 ordinal or dichotomous covariates. Again, the approximations based on the normal distribution were found to hold well, despite the clear deviations from the underlying assumptions.

Later work by Hill, Rubin and Thomas [8] also showed that the Rubin and Thomas [18] approximations held quite well with real data in the context of an evaluation of the New York School Choice Scholarship Program, which utilized randomization to award scholarships to eligible participants. Out of the large pool of possible controls, a matched sample was chosen for follow-up, where the matching was done using an affinely invariant matching method based on 21 ordinal or dichotomous covariates. Hill, Rubin and Thomas compared the bias and variance benefits of choosing matched controls rather than a random sample of controls. The Rubin and Thomas [18] results predict a gain of efficiency for differences in covariate means by a factor of approximately two, and Hill et al. showed that this predicted gain in efficiency was achieved, despite the markedly nonnormal distributions of some of the covariates.

In this paper we generalize the results of Rubin and Thomas [16, 17, 18] to the setting where the treated and control groups' covariate distributions are "discriminant mixtures of proportional ellipsoidally symmetric" (DMPES) distributions. We see that most, but not all, of the basic results in fact hold under these more general conditions, which support the broader applicability of these results, as suggested by the empirical evidence referenced above. We use as a running example the estimation of the effects of smoking on lung cancer, where the results here were used to motivate diagnostics for the results of matching [15].

**2. Discriminant mixtures of ellipsoidally symmetric distributions.** An ellipsoidal distribution for $p$-component $X$ is a distribution such that a linear transformation of $X$ leads to a spherically symmetric distribution, which is defined by the distribution on the radii of concentric hyperspheres on which there is a uniform probability density. Thus, an ellipsoidal distribution is specified by its center, inner product and distribution on the radius [5].

DEFINITION. The distribution on $X$, $F(X)$, is a "discriminant mixture of proportional ellipsoidally symmetric" (DMPES) distribution if it possesses the following properties:

(i) $F(X)$ is a mixture of $K$ ellipsoidally symmetric distributions $\{F_k; k = 1, \ldots, K\}$,

$$F(X) = \sum_{k=1}^{K} \alpha_k F_k(X), \tag{1}$$

where $\alpha_k \geq 0$ for all $k = 1, \ldots, K$, and $\sum_{k=1}^{K} \alpha_k = 1$, where $F_k$ has center $\mu_k$ and inner product $\Sigma_k$. Hence, the "mixture" (M) and "ellipsoidally symmetric" (ES) parts of DMPES.

(ii) The $K$ inner products are proportional:

$$\Sigma_i \propto \Sigma_j \quad \text{for all } i, j = 1, \ldots, K. \tag{2}$$

Hence, the "proportional" (P) part of DMPES.

(iii) The $K$ centers are such that all best linear discriminants between any two components are proportional:

$$(\mu_i - \mu_j)\Sigma_k^{-1} \propto (\mu_{i'} - \mu_{j'})\Sigma_{k'}^{-1} \quad \text{for all } i, j, k, i', j', k' = 1, \ldots, K. \tag{3}$$

Hence, the "discriminant" (D) in DMPES, because all mixture component centers lie along the common best linear discriminant.

In [16, 17, 18], $K = 2$, corresponding to the treated and control groups, and (2) is assumed; (3) is superfluous in the case with $K = 2$.

With DMPES distributions, there exists an affine transformation to a special canonical form, which is a simple extension of results in [3, 6] and [14]. This canonical form has, for each mixture component, the property that the distribution of $X$ is spherical, so that all inner products can be written as $\sigma_k^2 I$, where $I$ is the $p \times p$ identity matrix and $\sigma_k^2$ is a positive scalar constant, $k = 1, \ldots, K$. Moreover, the canonical form has the component centers lying along the unit vector (unless all $\mu_i = \mu_j$) so that the centers are $\delta_k U$, where $U = (1, \ldots, 1)'$, the $p$-component unit vector, and the $\delta_k$ are scalar constants, $k = 1, \ldots, K$; if all $\mu_i = \mu_j$, then all $\delta_k = 0$. Therefore, in their canonical form, the distribution of each component of $X$ is the same, and thus, the distribution of $X$ is exchangeable, not only within each of the $K$ mixture components, but also for any collection of mixture components defined by a subset of the indices $\{1, \ldots, K\}$.

Moreover, further symmetry results can be stated for a DMPES distribution by decomposing $X$ into its projection along the best linear discriminant, $Z$, and its projection orthogonal to $Z$. Specifically, the standardized best linear discriminant can be written as

$$Z = U'X/p^{1/2}, \tag{4}$$

unless all $\delta_k = 0$, in which case $Z$ is defined to be 0, the zero vector. Also, let $W$ be a standardized one-dimensional linear combination of $X$ orthogonal to $Z$,

$$(5) \qquad W = \gamma'X, \qquad \gamma'Z = 0, \qquad \gamma'\gamma = 1.$$

All such $W$ have the identical distribution in each mixture component, and the identical distribution for any collection of mixture components defined by a subset of the indices $\{1, \ldots, K\}$. Thus, the distribution of $X$ orthogonal to $Z$ has rotational symmetry, that is, is spherically symmetric.

Now suppose $K_t$ of the $K$ mixture components comprise the treatment group, and $K_c$ components comprise the control group, $K_t + K_c = K$; $K_t, K_c \geq 1$. Denote the set of treatment group component indices by $\mathcal{T}$ and the set of control group component indices by $\mathcal{C}$, $\mathcal{T} \cup \mathcal{C} = \{1, \ldots, K\}$. For example, $\mathcal{T}$ identifies current smokers and $\mathcal{C}$ identifies never smokers. The previous discussion implies that the distribution of $X$ is exchangeable in the treated group and in the control group, and moreover, the distribution of $X$ orthogonal to the discriminant $Z$ is spherically symmetric in the treated group and in the control group. This is the theoretical distributional setting for our results. In the more restrictive setting of [16] with proportional ellipsoidally symmetric distributions, $X$ is spherically symmetric in both groups.

**3. Results of matching with affinely invariant methods.** When affinely invariant matching methods are used with DMPES distributions, the canonical form given in Section 2 can be assumed without loss of generality. The following results, stated in canonical form, closely parallel results from [16]. The main symmetry arguments do not change with the use of mixtures of distributions. Although most of our results can be written without assuming finite first two moments in each mixture component and without restricting $K$ to be finite, the extra generality complicates notation and appears to be of little practical importance.

THEOREM 3.1. *Suppose an affinely invariant matching method is applied to random treated and control samples with DMPES distributions. Then*

$$E(\overline{X}_{mt}) \propto E(\overline{X}_{mc}) \propto U$$

*and*

$$\operatorname{var}(\overline{X}_{mt} - \overline{X}_{mc}) \propto I + cUU', \qquad c \geq -1/p,$$

*where $\overline{X}_{mt}$ and $\overline{X}_{mc}$ are the mean vectors in the matched treated and control samples, and $E(\cdot)$ and $\operatorname{var}(\cdot)$ are the expectation and variance over repeated random draws from the initial treated and control populations. Also,*

$$E(\nu_{mt}(X)) \propto I + c_t UU', \qquad c_t \geq -1/p,$$
$$E(\nu_{mc}(X)) \propto I + c_c UU', \qquad c_c \geq -1/p,$$

where $\nu_{mt}(X)$ and $\nu_{mc}(X)$ are the sample covariance matrices of $X$ in the matched treated and control groups, respectively. Corresponding formulas also hold within each of the mixture components. When $Z = 0$, $E(\overline{X}_{mt}) = E(\overline{X}_{mc}) = 0$, the zero vector, and $c = c_t = c_c = 0$.

PROOF. The proof follows directly from symmetry arguments and is essentially the same as that of Theorem 3.1 in [16]. Briefly, with affinely invariant matching methods, the matching treats each coordinate of $X$ the same and, hence, the exchangeability of the DMPES distributions of $X$ in matched treated and control samples is not affected. Thus, the expectations of the matched sample means of all coordinates of $X$ must be the same and, hence, the expectation of $X$ must be proportional to $U$ in each matched group. Analogously, the covariance matrices of $X$ must be exchangeable in each matched group. The general form for the covariance matrix of exchangeable variables is proportional to $I + cUU'$, $c \geq -1/p$. When $Z = 0$, the direction $U$ is no different from any other, that is, there is complete rotational symmetry and, hence, the simplification. □

COROLLARY 3.1. *The quantities* $\operatorname{var}(\overline{W}_{mt} - \overline{W}_{mc})$, $E(\nu_{mt}(W))$ *and* $E(\nu_{mc}(W))$ *take the same three values for all standardized $W$ orthogonal to $Z$. In addition, for each mixture component, $E(\nu_{mk}(W))$ takes the same value for all $W$, where $\nu_{mk}(W)$ is the sample variance of $W$ in the matched mixture component $k \in \mathcal{T}$ or $\mathcal{C}$.*

PROOF. The corollary follows from the fact that, due to the rotational symmetry in matched samples implied by Theorem 3.1 orthogonal to the discriminant, any $W$ will have the same distribution. □

**4. The effects on a linear combination of $X$ of affinely invariant matching relative to random sampling.** As in [16, 17, 18, 19], it is natural to describe the results of matching by its effects on a linear combination of $X$, $Y = \beta'X$, where, for convenience, we assume $Y$ is standardized, $\beta'\beta = 1$. Any such $Y$ can be expressed as the sum of projections along and orthogonal to the best linear discriminant,

$$(6) \qquad Y = \rho Z + (1 - \rho^2)^{1/2} W,$$

where $\rho$ is the correlation between $Y$ and $Z$. When $Z = 0$, $Y = W$ and $\rho \equiv 0$.

It is also natural, as in [16, 17, 18, 19], to compare the results of the matching to random sampling done in an affinely invariant way, such as randomly sampling from the original treated and control groups, thereby sampling from each component in proportion to its fraction in the population [the $\alpha$'s in (1)], or randomly sampling from each component with fixed proportions, where the same fixed proportions would be used in matching. We will refer to the treated and control samples generated by any such random sampling by indices $rt$ and $rc$, respectively, where $N_{rt} = N_{mt}$ and $N_{rc} = N_{mc}$, but generally, of course, $N_t \geq N_{rt}$ and $N_c \geq N_{rc}$.

The following corollaries decompose the effects on $Y$ of affinely invariant matching on $X$ into the effects of the matching on $Z$ and on $W$, relative to random sampling. Assuming the formulation from Section 2, we have the following results.

COROLLARY 4.1. (a) *When* $E(\overline{Z}_{rt} - \overline{Z}_{rc}) \neq 0$, *the matching is equal percent bias reducing (EPBR), as defined by* [14],

$$\text{(7)} \qquad \frac{E(\overline{Y}_{mt} - \overline{Y}_{mc})}{E(\overline{Y}_{rt} - \overline{Y}_{rc})} = \frac{E(\overline{Z}_{mt} - \overline{Z}_{mc})}{E(\overline{Z}_{rt} - \overline{Z}_{rc})}.$$

*Because the right-hand side of the above equation takes the same value for all $Y$, the percent bias reduction is the same for all $Y$.*

(b) *When $Z = 0$, the numerator and denominator of both ratios in equation (7) are 0.*

(c) *When $E(\overline{Z}_{rt} - \overline{Z}_{rc}) = 0$ but $Z \neq 0$, the denominators of both ratios in equation (7) are 0, and then $E(\overline{Y}_{mt} - \overline{Y}_{mc}) = \rho E(\overline{Z}_{mt} - \overline{Z}_{mc})$.*

PROOF. The proof of result (a) parallels the proof of Corollary 3.2 in [16]; however, here, rather than simple averages of $Z$, $W$ and $Y$, the averages are weighted averages of the mixture components, weighted, for example, by the $\alpha$'s in (1). Using the definition of $Y$,

$$E(\overline{Y}_{mt} - \overline{Y}_{mc}) = \rho E(\overline{Z}_{mt} - \overline{Z}_{mc}) + (\sqrt{1 - \rho^2}) E(\overline{W}_{mt} - \overline{W}_{mc}),$$

where, by the definition of $W$, $E(\overline{W}_{mt} - \overline{W}_{mc}) = \gamma' E(\overline{X}_{mt} - \overline{X}_{mc})$. From Theorem 3.1, $E(\overline{X}_{mt} - \overline{X}_{mc}) \propto U$ and again from the definition of $W$ in equation (5), $\gamma'Z = 0$. Thus,

$$E(\overline{Y}_{mt} - \overline{Y}_{mc}) = \rho E(\overline{Z}_{mt} - \overline{Z}_{mc}).$$

Similarly, $E(\overline{Y}_{rt} - \overline{Y}_{rc}) = \rho E(\overline{Z}_{rt} - \overline{Z}_{rc})$ because $E(\overline{W}_{rt} - \overline{W}_{rc}) = 0$ and result (a) of Corollary 4.1 follows.

Results (b) and (c) follow by analogous arguments. Situation (c) cannot arise when $K = 2$ because, with only one treated and one control component, $E(\overline{Z}_{rt} - \overline{Z}_{rc}) = 0$ implies that $Z = 0$. However, with multiple components in the treated and control groups, the difference in weighted averages ($E(\overline{Z}_{rt} - \overline{Z}_{rc})$) can equal 0 without all of the mixture component centers ($\{\mu_k\}$) being 0. □

This corollary implies that affinely invariant matching that reduces bias in one direction cannot create bias in some other direction. If bias reduction is obtained along $Z$, it is also obtained for all $Y$.

COROLLARY 4.2. *The matching is $\rho^2$ proportionate modifying of the variance of the difference in matched sample means,*

$$(8) \quad \frac{\operatorname{var}(\overline{Y}_{mt} - \overline{Y}_{mc})}{\operatorname{var}(\overline{Y}_{rt} - \overline{Y}_{rc})} = \rho^2 \frac{\operatorname{var}(\overline{Z}_{mt} - \overline{Z}_{mc})}{\operatorname{var}(\overline{Z}_{rt} - \overline{Z}_{rc})} + (1 - \rho^2) \frac{\operatorname{var}(\overline{W}_{mt} - \overline{W}_{mc})}{\operatorname{var}(\overline{W}_{rt} - \overline{W}_{rc})},$$

*where the ratios*

$$\frac{\operatorname{var}(\overline{Z}_{mt} - \overline{Z}_{mc})}{\operatorname{var}(\overline{Z}_{rt} - \overline{Z}_{rc})} \quad and \quad \frac{\operatorname{var}(\overline{W}_{mt} - \overline{W}_{mc})}{\operatorname{var}(\overline{W}_{rt} - \overline{W}_{rc})}$$

*take the same two values for all $Y$.*

PROOF. Using the definitions of $Z$ and $W$ in (4) and (5),

$$\operatorname{cov}(\overline{Z}_{mt} - \overline{Z}_{mc}, \overline{W}_{mt} - \overline{W}_{mc}) = \frac{1}{\sqrt{p}} U' \operatorname{var}(\overline{X}_{mt} - \overline{X}_{mc}) \gamma,$$

which from Theorem 3.1 is proportional to

$$U'(I + cUU')\gamma = U'\gamma + cpU'\gamma = 0,$$

again using the definition of $W$ in (5). Then, from the definition of $Y$ in (6),

$$\operatorname{var}(\overline{Y}_{mt} - \overline{Y}_{mc}) = \rho^2 \operatorname{var}(\overline{Z}_{mt} - \overline{Z}_{mc}) + (1 - \rho^2) \operatorname{var}(\overline{W}_{mt} - \overline{W}_{mc}).$$

Equation (8) follows because, in random subsamples, the samples from each treated and control mixture component are independent with

$$\operatorname{var}(\overline{Y}_{rt} - \overline{Y}_{rc}) = \operatorname{var}(\overline{Z}_{rt} - \overline{Z}_{rc}) = \operatorname{var}(\overline{W}_{rt} - \overline{W}_{rc}).$$

Also, $\operatorname{var}(\overline{Y}_{rt}) = \operatorname{var}(\overline{Z}_{rt})$ and $\operatorname{var}(\overline{Y}_{mc}) = \operatorname{var}(\overline{Z}_{mc})$, and each is a weighted linear combination of the variances in each of the treated and control mixture components, respectively. The final statement of Corollary 4.2 follows from Corollary 3.1. □

COROLLARY 4.3. *Within each of the mixture components, the matching is $\rho^2$ proportionate modifying of the expectation of the sample variances,*

$$(9) \quad \frac{E(\nu_{mk}(Y))}{E(\nu_{rk}(Y))} = \rho^2 \frac{E(\nu_{mk}(Z))}{E(\nu_{rk}(Z))} + (1 - \rho^2) \frac{E(\nu_{mk}(W))}{E(\nu_{rk}(W))},$$

*where $\nu_{rk}(\cdot)$ is the sample variance of $n_k$ randomly chosen units from component $k$, and $\nu_{mk}(\cdot)$ is the sample variance of $n_k$ matched units from component $k$ ($k \in \mathcal{T}$ or $\mathcal{C}$), and the ratio*

$$\frac{E(\nu_{mk}(W))}{E(\nu_{rk}(W))}$$

*takes the same value for all $Y$. The same is true for $E(\nu_{mk}(Z))/E(\nu_{rk}(Z))$.*

PROOF. In the matched sample from component $k \in \mathcal{T}$ or $\mathcal{C}$, the expected covariance of $Z$ and $W$ is

$$E(\text{cov}_{mk}(Z,W)) = \frac{1}{\sqrt{p}} E(U' \nu_{mk}(X) \gamma) \propto U'(I + c_k p U U') \gamma = 0,$$

from Theorem 3.1 and the definition of $W$ in (5), and $\nu_{mk}(X) \propto I + c_k UU'$, where the constants $c_k \geq -1/p$. Then, from (6),

$$E(\nu_{mk}(Y)) = \rho^2 E(\nu_{mk}(Z)) + (1-\rho^2) E(\nu_{mk}(W)).$$

Equation (9) follows because $E(\nu_{rk}(Y)) = E(\nu_{rk}(Z)) = E(\nu_{rk}(W))$. The final statement follows from Corollary 3.1. □

Note that the version of Corollary 4.3 stated for the full treated and control groups does not hold. In the special case considered in [16], there is only one component in each group.

**5. Conditionally affinely invariant matching with conditionally DMPES distributions.** We now extend the results of the previous sections to a setting where a subset of the covariates is treated differently from the remainder of the covariates, for example, exact matching on gender followed by discriminant matching, or Mahalanobis matching on key covariates within propensity score calipers [13]. Such matching was done, for example, in [15] when creating matched samples of current smokers and never smokers.

We define $X^{(s)}$ to be the $s$ "special covariates" spanning an $s$-dimensional subspace (e.g., gender, race in the smoking example) and $X^{(r)}$ to be the $r = p - s$ remaining covariates spanning an $r$-dimensional subspace (e.g., education, age). The methods considered are "conditionally affinely invariant matching methods" [16], which have the property that the result of the matching is the same following any (full-rank) affine transformation of the "remainder" covariates $X^{(r)}$:

$$((\mathcal{X}_t^{(s)}, \mathcal{X}_t^{(r)}), (\mathcal{X}_c^{(s)}, \mathcal{X}_c^{(r)})) \to (T, C)$$

implies

$$((\mathcal{X}_t^{(s)}, A(\mathcal{X}_t^{(r)})), (\mathcal{X}_c^{(s)}, A(\mathcal{X}_c^{(r)}))) \to (T, C).$$

In parallel with Section 2, we consider the case where each mixture component of the full covariate distribution has mean vectors $\mu_k^{(s)}$ and $\mu_k^{(r)}$, covariance matrices $\Sigma_k^{(s)}$ and $\Sigma_k^{(r)}$, and conditional means and covariance matrices given by $\mu_k^{(r|s)}$ and $\Sigma_k^{(r|s)}$. The full distribution of $X = (X^{(r)}, X^{(s)})$ across both groups is a conditional DMPES distribution if (i) the conditional distribution $X^{(r)}|X^{(s)}$ is ellipsoidal in each mixture component, (ii) it has proportional conditional covariance matrices, $\Sigma_k^{(r|s)} \propto \Sigma_{k'}^{(r|s)}$ for all $k$ and $k'$, and

(iii) it has centers such that $(\mu_i^{(r|s)} - \mu_j^{(r|s)})\Sigma_k^{(r|s)-1} \propto (\mu_{i'}^{(r|s)} - \mu_{j'}^{(r|s)})\Sigma_{k'}^{(r|s)-1}$ for all $i, j, k, i', j', k' = 1, \ldots, K$. Notice that condition (ii) implies a common (across all mixture components) linear regression of the $r$ covariates in $X^{(r)}$ on the $s$ covariates in $X^{(s)}$, with coefficients $B$. As noted by [16], the special case with $X^{(s)}$ binomial or multinomial and $X^{(r)}$ multivariate normal relates to the logistic regression model for predicting treated or control status given the covariates, thus relating it to the methods of propensity score estimation developed by Rosenbaum and Rubin [12, 13].

We again can use a canonical form when a conditionally affinely invariant matching method is used with a conditionally DMPES distribution. The covariates $X^{(r)}$ are redefined as the components of $X^{(r)}$ uncorrelated with $X^{(s)}$: $X^{(r)} - B'X^{(s)}$. The following notation is then used for the moments of the distribution of $X^{(r)}$ [and the conditional moments of $X^{(r)}$ given $X^{(s)}$]:

$$\mu_k^{(r)} = \delta_k^{(r)} U, \qquad \Sigma_k^{(r)} = \sigma_k^2 I,$$

$k = 1, \ldots, K$, where $\delta_k^{(r)}$ and $\sigma_k^2$ are scalar constants, $U$ is now the $r$-dimensional unit vector, and $I$ is now the $r \times r$ identity matrix. Thus, the distributions of $(X^{(s)}, X^{(r)})$ and $X^{(r)}$ given $X^{(s)}$ are exchangeable under permutations of components of $X^{(r)}$ conditional on $X^{(s)}$ in each of the mixture components.

THEOREM 5.1. *Suppose a conditionally affinely invariant matching method is applied to random treated and control samples with conditional DMPES distributions. Then, in canonical form,*

$$E(\overline{X}_{mt}^{(r)}) \propto U, \qquad E(\overline{X}_{mc}^{(r)}) \propto U$$

*and*

$$\operatorname{var}(\overline{X}_{mt} - \overline{X}_{mc}) = \begin{bmatrix} \operatorname{var}(\overline{X}_{mt}^{(s)} - \overline{X}_{mc}^{(s)}) & CU' \\ UC' & k(I + c_0 UU') \end{bmatrix},$$

*where $k \geq 0$, $c_0 \geq -1/r$ and $C' = (c_1, \ldots, c_s)$. Also,*

$$E(\nu_{mt}(X)) = \begin{bmatrix} E(\nu_{mt}(X^{(s)})) & C_t U' \\ UC_t' & k_t(I + c_{t0} UU') \end{bmatrix},$$

*where $k_t \geq 0$, $c_{t0} \geq -1/r$, $C_t' = (c_{t1}, c_{t2}, \ldots, c_{ts})$, with an analogous result and notation for the matched control group. When $Z = 0$, $E(\overline{X}_{mt}^{(r)}) = E(\overline{X}_{mc}^{(r)}) = 0$, $C = C_t = C_c = 0$, the zero vector, and $c_0 = c_{t0} = c_{c0} = 0$.*

PROOF. The proof of this theorem parallels that of Theorem 3.1, with the exception of the existence of the covariances between components in $X^{(s)}$ and $X^{(r)}$. Due to the symmetry, these covariances are also exchangeable in the coordinates of $X^{(r)}$. □

**6. Effect on $Y$ of matching with special covariates.** In parallel with the earlier formulation, we express an arbitrary linear combination of $X$ as

$$Y = \rho \mathcal{Z} + (1-\rho^2)^{1/2}\mathcal{W},$$

where $\mathcal{Z}$ and $\mathcal{W}$ are the standardized projections of $Y$ along and orthogonal to the subspace spanned by $(X^{(s)}, Z)$, respectively, and $\rho$ is the correlation between $Y$ and $\mathcal{Z}$. In this framework, $Z$ is the standardized discriminant of the covariates uncorrelated with $X^{(s)}$, again expressed in canonical form as $Z = U'X^{(r)}/r^{1/2}$. When $\mu_k^{(r)} = 0$ for all $k$, $Z$ is defined to be the zero vector, and then $\mathcal{Z}$ is defined to be the projection of $Y$ in the subspace spanned by $X^{(s)}$.

We write $\mathcal{Z}$ and $\mathcal{W}$ as

$$\mathcal{Z} = \psi'X = (\psi^{(s)\prime}, \psi^{(r)\prime})\begin{pmatrix} X^{(s)} \\ X^{(r)} \end{pmatrix}, \tag{10}$$

$$\mathcal{W} = \gamma'X = (\gamma^{(s)\prime}, \gamma^{(r)\prime})\begin{pmatrix} X^{(s)} \\ X^{(r)} \end{pmatrix}. \tag{11}$$

LEMMA 6.1. *The coefficients $\gamma$ and $\psi$ satisfy*

$$\gamma^{(s)} = (0,\ldots,0)', \qquad \gamma^{(r)\prime}\psi^{(r)} = Z\gamma^{(r)} = 0 \quad \text{and} \quad \psi^{(r)} \propto U. \tag{12}$$

PROOF. Equation (12) follows because $\mathcal{W}$ is a linear combination of $X$ uncorrelated with $\{X^{(s)}, Z\}$, and thus uncorrelated with $\{X^{(s)}\}$, and because $\mathcal{Z}$ is uncorrelated with $\mathcal{W}$. The other results follow from these and the definition of $\mathcal{Z}$ in canonical form. □

Because the symmetry results of Theorem 5.1 for $X$ orthogonal to $\mathcal{Z}$ imply that all $\mathcal{W}$ orthogonal to $\mathcal{Z}$ have the same distribution, we immediately have the following corollary to Theorem 5.1.

COROLLARY 6.1. *The quantities* $\text{var}(\mathcal{W}_{mt} - \mathcal{W}_{mc})$, $E(\nu_{mt}(\mathcal{W}))$, *and* $E(\nu_{mc}(\mathcal{W}))$ *take the same three values for all standardized* $Y$. *Analogous results hold for statistics in random subsamples indexed by $rt$ and $rc$. In addition, $E(\nu_{mk}(\mathcal{W}))$ takes the same value for all $\mathcal{W}$ within each of the mixture components, $k \in \mathcal{T}$ or $\mathcal{C}$. However, the corresponding expressions involving $\mathcal{Z}$ generally do depend on the choice of $Y$.*

COROLLARY 6.2. (a) *When $E(\overline{\mathcal{Z}}_{rt} - \overline{\mathcal{Z}}_{rc}) \neq 0$, the percent bias reduction in $Y$ equals the percent bias reduction of $Y$ in the subspace $\{X^{(s)}, Z\}$,*

$$\frac{E(\overline{Y}_{mt} - \overline{Y}_{mc})}{E(\overline{Y}_{rt} - \overline{Y}_{rc})} = \frac{E(\overline{\mathcal{Z}}_{mt} - \overline{\mathcal{Z}}_{mc})}{E(\overline{\mathcal{Z}}_{rt} - \overline{\mathcal{Z}}_{rc})}.$$

(b) *When $E(\overline{\mathcal{Z}}_{rt} - \overline{\mathcal{Z}}_{rc}) = 0$, the denominators of both ratios in* (a) *equal* 0, *and* $E(\overline{Y}_{mt} - \overline{Y}_{mc}) = \rho E(\overline{\mathcal{Z}}_{mt} - \overline{\mathcal{Z}}_{mc})$.

PROOF. The proof parallels that of Corollary 4.1 because $\mathcal{W} = \gamma'^{(r)} X^{(r)}$ from the definition of $\mathcal{W}$ in (11) and Lemma 6.1, and from Theorem 5.1 and Lemma 6.1, $\gamma'^{(r)} E(\overline{X}_{mt}^{(r)} - \overline{X}_{mc}^{(r)}) = 0$. Thus, $E(\overline{\mathcal{W}}_{mt} - \overline{\mathcal{W}}_{mc}) = 0$.  □

COROLLARY 6.3. *The matching is $\rho^2$ proportionate modifying of the variance of the difference in matched sample means,*

$$\frac{\operatorname{var}(\overline{Y}_{mt} - \overline{Y}_{mc})}{\operatorname{var}(\overline{Y}_{rt} - \overline{Y}_{rc})} = \rho^2 \frac{\operatorname{var}(\overline{\mathcal{Z}}_{mt} - \overline{\mathcal{Z}}_{mc})}{\operatorname{var}(\overline{\mathcal{Z}}_{rt} - \overline{\mathcal{Z}}_{rc})} + (1 - \rho^2) \frac{\operatorname{var}(\overline{\mathcal{W}}_{mt} - \overline{\mathcal{W}}_{mc})}{\operatorname{var}(\overline{\mathcal{W}}_{rt} - \overline{\mathcal{W}}_{rc})},$$

*where the ratio* $\operatorname{var}(\overline{\mathcal{W}}_{mt} - \overline{\mathcal{W}}_{mc})/\operatorname{var}(\overline{\mathcal{W}}_{rt} - \overline{\mathcal{W}}_{rc})$ *takes the same value for all* $Y$.

PROOF. The proof is analogous to that of Corollary 4.2 using Theorem 5.1 and Lemma 6.1, and parallels the proof of Corollary 4.3 in [16], where, in that proof, there is a typographical error: $\overline{Z}_{mt} - \overline{Z}_{mc}$ and $\overline{W}_{mt} - \overline{W}_{mc}$ should be replaced by $\overline{\mathcal{Z}}_{mt} - \overline{\mathcal{Z}}_{mc}$ and $\overline{\mathcal{W}}_{mt} - \overline{\mathcal{W}}_{mc}$, respectively.  □

COROLLARY 6.4. *Within each mixture component, the matching is $\rho^2$ proportionate modifying of the expectation of the sample variances,*

$$\frac{E(\nu_{mk}(Y))}{E(\nu_{rk}(Y))} = \rho^2 \frac{E(\nu_{mk}(\mathcal{Z}))}{E(\nu_{rk}(\mathcal{Z}))} + (1 - \rho^2) \frac{E(\nu_{mk}(\mathcal{W}))}{E(\nu_{rk}(\mathcal{W}))}$$

*for all* $k \in \mathcal{T}$ *or* $\mathcal{C}$, *where the ratio* $E(\nu_{mk}(\mathcal{W}))/E(\nu_{rk}(\mathcal{W}))$ *takes the same value for all* $Y$ *within each mixture component.*

PROOF. The proof of this corollary parallels that of Corollary 4.3, with modifications similar to those in the proof of Corollary 6.3. Again, as in Corollary 4.2, this result generally holds only in each of the individual treated and control group components, and the analogous result in the overall samples does not hold.  □

**7. Discussion.** Here we have shown that most of the results proven by Rubin and Thomas [16] can be extended to discriminant mixtures of proportional ellipsoidally symmetric (DMPES) distributions, as defined in Section 2, and provides some theoretical rationale for why the earlier Rubin and Thomas [16, 17, 18] results hold well even when the assumption of ellipsoidally symmetric distributions is not met. These results show that even with the more complicated setting of DMPES distributions, the effects of matching on an arbitrary linear combination of the covariates can be summarized by its effects along and orthogonal to the discriminant.

Although the class of DMPES distributions is still restrictive, previous experience has indicated that mathematically convenient conditions for matching can provide guidance in real-world examples. A classic example is in [2] on the bias reduction possible from stratified matching. Although Cochran's results were proved assuming infinite samples sizes and a linear relationship between a single covariate and the outcome, the approximations and their implied guidance have found applicability and use for a much wider range of situations. For a specific example here, the implications of our results were the basis for the applied diagnostics in [15] used to assess the quality of the matched samples of smokers and never smokers in the National Medical Expenditure Survey, based on decomposing the comparisons of the distributions in the matched samples into components along and orthogonal to the discriminant.

**Acknowledgments.** The authors would like to thank the previous and current Editors, anonymous referees and Neal Thomas for exceptionally helpful comments. This work was done while Elizabeth Stuart was a graduate student in the Department of Statistics, Harvard University.

DEPARTMENT OF STATISTICS
HARVARD UNIVERSITY
1 OXFORD ST.
CAMBRIDGE, MASSACHUSETTS 02138
USA
E-MAIL: rubin@stat.harvard.edu

JOHNS HOPKINS BLOOMBERG
SCHOOL OF PUBLIC HEALTH
624 N. BROADWAY, ROOM 804
BALTIMORE, MARYLAND 21205
USA
E-MAIL: estuart@jhsph.edu